\documentclass[12pt,reqno]{article}
\usepackage[english]{babel}
\usepackage{amsmath,amsthm}
\usepackage{amsfonts}
\usepackage{graphicx}
\usepackage{enumerate}
\usepackage[usenames,dvipsnames]{color}
\usepackage{subfigure}
\usepackage[right,pagewise,displaymath, mathlines]{lineno}
\usepackage{epstopdf}
\usepackage{color}
\usepackage{multirow}

\numberwithin{equation}{section}
\numberwithin{figure}{section}
\numberwithin{table}{section}

\newtheorem{theorem}{Theorem}[section]
\newtheorem{corollary}{Corollary}[section]

\numberwithin{equation}{section}

\begin{document}

\begin{center}

{\Large\bf Weighted allocations, their concomitant-based estimators, and asymptotics}

\vspace*{7mm}

{\large Nadezhda Gribkova}

\medskip

\textit{Faculty of Mathematics and Mechanics, St.\,Petersburg State University,\\ St.\,Petersburg 199034, Russia}

\bigskip

{\large Ri\v cardas Zitikis}

\medskip

\textit{School of Mathematical and Statistical Sciences,
Western University, \break London, Ontario N6A 5B7, Canada}

\end{center}

\begin{abstract}
Various members of the class of weighted insurance premiums and risk capital allocation rules have been researched from a number of perspectives. Corresponding formulas in the case of parametric families of distributions have been derived, and they have played a pivotal role when establishing parametric statistical inference in the area. Non-parametric inference results have also been derived in special cases such as the tail conditional expectation, distortion risk measure, and several members of the class of weighted premiums. For weighted allocation rules, however, non-parametric inference results have not yet been adequately developed. In the present paper, therefore, we put forward empirical estimators for the weighted allocation rules and establish their consistency and asymptotic normality under practically sound conditions. Intricate statistical considerations rely on the theory of induced order statistics, known as concomitants.

\medskip

\noindent
{\it Key words and phrases:} weighted allocation, premium, consistency, asymptotic normality, concomitant.

\medskip

\noindent
{\it 2010 MSC:}  Primary: 62G05, 62G20;
Secondary: 62P05, 62P20.
\end{abstract}

\section{Introduction}
\label{intro}

The tail conditional expectation, exponential tilting, and various members of the class of weighted insurance premiums and the corresponding risk capital allocation rules have been extensively researched (e.g., Pflug and R\"{o}misch, 2007; R\"{u}schendorf, 2013; McNeil et al., 2015; F\"{o}llmer and Schied, 2016). Their formulas in the case of various parametric families of distributions have been derived (e.g., Furman and Landsman, 2005, 2010; Asimit et al., 2013; Su, 2016; Asimit et al., 2016; Su and Furman, 2017; Ratovomirija et al., 2017; Vernic, 2017; and references therein), thus facilitating parametric statistical inference in the area. The literature also contains a number of non-parametric inference results (e.g., Brazauskas and Serfling, 2003; Brazauskas, 2009; Brazauskas and Kleefeld, 2009; and references therein), particularly in special cases such as the value at risk (e.g., Maesono and Penev, 2013; and references therein), the tail conditional expectation (e.g., Brazauskas et al., 2008; and references therein), and distortion risk measures (e.g., Jones and Zitikis, 2007; and references therein), with results available in light- and heavy-tailed settings (e.g., Necir and Meraghni,  2009; Necir et al., 2007; Necir et al., 2010; Rassoul, 2013; Brahimi et al., 2012; and references therein). Non-parametric statistical inference for weighted allocation rules has not yet, however, been adequately developed, and we therefore devote the current paper to this topic. In particular, we construct empirical estimators for the weighted allocations and establish their consistency and asymptotic normality under practically sound conditions. Details follow.

Let $X$ be a real-valued random variable, which could, for example, be a financial or insurance risk associate with a business line in a company. Denote the cumulative distribution function (cdf) of $X$ by $F_X$. When $X$ is viewed as a stand-alone risk, then the capital needed to mitigate the risk can be calculated using (e.g., Furman and Zitikis, 2008a)
\begin{equation}
\label{0pi_w}
\pi_w=\frac{\textbf{E}[X w\circ F_X(X)]}{\textbf{E}[w\circ F_X(X)]}
\end{equation}
with an appropriately chosen weight function $w:[0,1]\to [-\infty,\infty ]$, where  $w\circ F_X$ denotes the composition of the functions $w$ and $F_X$. We of course assume that the two expectation in the definition of $\pi_w$ are well-defined and finite, and $\textbf{E}[w\circ F_X(X)]$ is not zero.

The function $w$ may or may not take infinite value, may or may not be non-decreasing, depending on the context. Throughout the paper we always assume that $w$ is finite on the open interval $(0,1)$ and, at each point of $(0,1)$, is either left-continuous or right-continuous. As far as we are aware of, all practically relevant weight functions satisfy these properties, with a few illustrative examples following next.

When dealing with insurance losses, researchers work with non-negative and non-decreasing weight functions, which ensure that $\pi_w$ is non-negatively loaded, that is, the bound $\pi_w\ge \textbf{E}[X]$ holds for all risks $X$ under consideration. In other contexts, such as econometrics and, more specifically, measurement of income inequality, the function $w$ can be non-increasing.

For example, $w(t)=\mathbf{1}\{t>p \}$ for any parameter $p\in (0,1)$ is non-decreasing and leads to the insurance version of the tail conditional expectation. Another example is $w(t)=\nu(1-t)^{\nu-1}$ with  parameter $\nu>0$. If $\nu \in (0,1]$, then $\pi_w$ is the proportional hazards transform (Wang, 1995, 1996). If $\nu \ge 1$, then $\pi_w$ reduces to the (absolute) $S$-Gini index used for measuring income equality (e.g., Zitikis and Gastwirth, 2002; and references therein).

Note that if the cdf $F_X$ is continuous, then $\pi_w$ can be written as the integral
\begin{equation}
\label{0pi_w-small}
\pi_w=\int_0^1 F_X^{-1}(t)w^*(t)\mathrm{d}t
\end{equation}
of the quantile function $F_X^{-1}$ of $X$, with the weight function
\[
w^*(t)={w(t)\over \int_0^1 w(u)\mathrm{d}u},
\]
which is a probability density function (pdf) whenever $w(t)\ge 0$ for all $t\in[0,1]$ and $\int_0^1 w(u)\mathrm{d}u\in (0,\infty )$.  This represetation of $\pi_w$ connects our present research with the dual utility theory (Yaari, 1987; Quiggin, 1993; and references therein) that has arisen as a prominent alternative to the classical utility theory of von Neumann and Morgenstern (1944).

Setting appropriate insurance premiums and allocating capital to individual business lines are usually done within a company's risk profile. That is, if $Y$ is the risk associated with the entire company, then allocating capital to the business line whose risk is $X$ is done by taking into account the value of $Y$. This viewpoint leads us to the weighted capital allocation rule  (Furman and Zitikis, 2008b)
\begin{equation}
\label{pi_w}
\Pi_w=\frac{\textbf{E}[X w\circ F_Y(Y)]}{\textbf{E}[w\circ F_Y(Y)]},
\end{equation}
where $F_Y$ denotes the cdf of $Y$. Obviously, setting $Y$ to $X$ reduces $\Pi_w$ to $\pi_w$, and for this reason we concentrate on developing  nonparametric statistical inference for $\Pi_w$ and then specialize our results to $\pi_w$. For the role of $\pi_w$ and $\Pi_w$ in the context of the weighted insurance pricing model (WIPM), we refer to Furman and Zitikis (2017).

To construct an empirical estimator for $\Pi_w$, let $(X_k,Y_k)$, $k=1,2\dots$, be independent copies of the random pair $(X,Y)$ and, for each integer $n\geq 1$, let $\widehat{F}_Y$ be defined by
\begin{equation}
\label{cdf-y}
\widehat{F}_Y(y)={1\over n+1}\sum_{k=1}^{n} \mathbf{1}\{Y_k\le y\},
\end{equation}
where $\mathbf{1}\{Y_k\le y\}$ is the indicator of $Y_k\le y$: it is equal to $1$ when $Y_k\le y$ is true, and $0$ otherwise. This is an empirical estimator of the cdf $F_Y(y)$ that slightly differs from the classical empirical cdf because we use $1/(n+1)$ instead of $1/n$. This adjustments is important as in this way defined $\widehat{F}_Y(y)$ takes only values $k/(n+1)$, $k=1, \dots , n$, which are always inside the open interval $(0,1)$ on which the weight function $w$ is finite.

We are now in the position to define the empirical estimator of  $\Pi_w$ by the formula
\begin{equation}
\label{plug_in_1}
\widetilde{\Pi}_{w}=\frac{\sum_{k=1}^{n}X_k w\circ \widehat{F}_Y(Y_k)}{\sum_{k=1}^{n} w\circ \widehat{F}_Y(Y_k)},
\end{equation}
with a tilde used instead of the usual hat on top of $\Pi_{w}$ because we reserve the latter notation for another estimator to be introduced in a moment. Note that when $Y=X$ and thus $Y_k=X_k$ for all $k\ge 1$, the empirical allocation $\widetilde{\Pi}_{w}$ reduces to the estimator
\begin{equation}
\label{plug_in_1b}
\widetilde{\pi}_{w}=\frac{\sum_{k=1}^{n}X_k w\circ \widehat{F}_X(X_k)}{\sum_{k=1}^{n} w\circ \widehat{F}_X(X_k)}
\end{equation}
of $\pi_w$, where $\widehat{F}_X$ is defined by equation (\ref{cdf-y}) but now based on $X_1, \dots, X_n$. Estimators (\ref{plug_in_1}) and (\ref{plug_in_1b}) are ratio statistics, whose deep asymptotic properties have been explored by Maesono (2005, 2010).

When the underlying population cdf $F_Y$ is continuous, the random variable $F_Y(Y)$ is uniform on $(0,1)$ and, therefore, the denominator in the definition of $\Pi_w$ is equal to $\int_0^1w(u)\,\mathrm{d}u$. Since we do not need to estimate the latter integral, we can therefore use the following simpler estimator
\begin{equation}
\label{plug_in_2}
\widehat{\Pi}_{w}={\widehat{\Delta}_{w}\over \int_0^1 w(u)\mathrm{d}u}
\end{equation}
of $\widehat{\Pi}_{w}$, where
\[
\widehat{\Delta}_{w}
={1\over n} \sum_{k=1}^{n}X_k w\circ \widehat{F}_Y(Y_k) .
\]
Note that
\[
\widehat{\Delta}_{w}
\stackrel{a.s.}{=}{1\over n} \sum_{k=1}^{n}X_{[k:n]} w_{k,n},
\]
where
\[
w_{k,n}=w\bigg ( {k\over n+1} \bigg )
\]
and $X_{[1:n]}, \dots , X_{[n:n]}$ are the induced order statistics, known as concomitants, corresponding to $Y_{1:n}, \dots , Y_{n:n}$. When $Y=X$ and thus $Y_k=X_k$ for all $k\ge 1$, then the concomitants reduce to the order statistics $X_{1:n}, \dots , X_{n:n}$. In this case, the estimator $\widehat{\Pi}_{w}$ reduces to the estimator $\widehat{\pi}_{w}$ of $\pi_w$  given by the equation
\begin{equation}
\label{plug_in_2-cor}
\widehat{\pi}_{w}={1\over n} \sum_{k=1}^{n}X_{k:n} {w_{k,n}\over \int_0^1 w(u)\mathrm{d}u}
\end{equation}

While $\widehat{\pi}_{w}$ is a linear combination of order statistics (e.g., Gribkova, 2017; and references therein), which is a less technically demanding object, the estimators $\widetilde{\Pi}_{w}$,
$\widehat{\Pi}_{w}$, and $\widehat{\Delta}_{w}$ are linear combinations of concomitants, which require much more sophisticated methods of analysis. In what follows, we establish conditions under which these estimators are consistent and asymptotically normal. Main results are in Section \ref{sect-2}, with their proofs in Section \ref{sect-3a}. Throughout the paper we use $c$ to denote various constants that do not depend on $n$ and usually change their values from line to line. Furthermore, we use $\stackrel{a.s.}{\to}$ to denote convergence almost surely, $\stackrel{\mathbf{P}}{\to}$ convergence in probability, and $\stackrel{law}{\to}$ convergence in law/distribution. We use ``:='' when  wishing to emphasize that a certain equation holds by definition.

\section{Main results}
\label{sect-2}

Three quantities influence the asymptotic behaviour of the above introduced estimators: 1) the weight function $w$, 2) the cdf of $X$, and 3) the cdf of $Y$. They interact with each other, and thus determining their influence on asymptotic results (e.g., consistency, asymptotic normality, etc.) in one stroke becomes not only challenging but also leads to unwieldy -- from the practical point of view -- conditions. Because of this issue, we next set out to establish asymptotic results in several complementary forms, starting with strong consistency.

\begin{theorem}
\label{theorem-1}
If the first moment $\mathbf{E}[X]$ is finite and the weight function $w$ is continuous on $[0,1]$, then
$\widetilde{\Pi}_{w} \stackrel{a.s.}{\to} \Pi_w$
and thus $\widetilde{\pi}_{w} \stackrel{a.s.}{\to} \pi_w $ when $n\to\infty$.
\end{theorem}

The theorem is attractive in the sense that it does not impose any condition on the underlying random variables, except the very minimal condition that the first moment of $X$ is finite. The condition on the weight function $w$ is, however, very strong. For example, it is not satisfied by $w(t)=\nu(1-t)^{\nu-1}$ for any $\nu \in (0,1)$. Furthermore, the condition is not satisfied by $w(t)=\mathbf{1}\{t>p \}$ for any $p\in (0,1)$, and we thus cannot use the theorem to deduce strong consistency of the tail conditional expectation.

In the next two theorems we no longer assume continuity (and thus boundedness) of the weight function $w$ on the compact interval $[0,1]$. Instead, we require finite higher moments of $X$, as well as the continuity of the cdf $F_Y$ when dealing with $\widehat{\Pi}_{w}$ and the continuity of the cdf $F_X$ when dealing with $\widehat{\pi}_{w}$.

We use $L_p$, $1\leq p\leq \infty$, to denote the space of all Borel measurable functions $h: [0,1]\to \mathbb{R}$ such that $\|h\|_p :=(\int_0^1 |f|^p\mathrm{d}\lambda)^{1/p}<\infty $ when $1\leq p<\infty$ and $\|h\|_{\infty}:=\mathrm{ess}\sup_{t\in [0,1]}|h(t)|<\infty $ when $p=\infty$, where $\lambda $ is the Lebesgue measure.

The following two theorems are consequences of the strong law of large numbers for $L$-statistics, proved in various levels of generality by van~Zwet (1980). For example, Theorem 2.1 and Corollary 2.1 of van~Zwet (1980) imply the following theorem.

\begin{theorem}
\label{theorem-zwet-1}
Let the cdf $F_X$ be continuous.  If $\mathbf{E}[|X|^{p}]<\infty$ and $w\in L_{q}$ for some $p,q \in [1, \infty ]$ such that $p^{-1}+q^{-1}=1$, then
$\widetilde{\pi}_{w} \stackrel{a.s.}{\to} \pi_w $ when $n\to\infty$.
\end{theorem}

The next theorem, which follows from Theorem 3.1 of van~Zwet (1980), allows us to use different values of $p$ and $q$ on different subintervals of $(0,1)$, thus enabling different growth rates of the  quantile function $ F_X^{-1}$ and the weight function $w$ near the endpoints of the interval $(0,1)$.  Following van~Zwet (1980), let $0=:a_0<a_1<\dots <a_j:=1$ be points dividing the interval $(0,1)$  into $j\ge 1$ subintervals, which we denote by
\[
A_i=(a_{i-1},a_i), \quad i=1,\dots ,j,
\]
whose $\varepsilon$-neighbourhoods within the interval $(0,1)$  are
\[
B_{i,\varepsilon}=(a_{i-1}-\varepsilon,a_i +\varepsilon)\cap (0,1).
\]

\begin{theorem}
\label{theorem-zwet-2}
Let the cdf $F_X$ be continuous, and let $p_i,\, q_i \in [1,\infty ]$ be such that $p_i^{-1}+q_i^{-1}=1$. If there is $\varepsilon>0$ such that  $F^{-1}_X\mathbf{1}_{B_{i,\varepsilon}}\in L_{p_i}$ and $w \mathbf{1}_{A_i}\in L_{q_i}$ for every $i=1,\dots ,j$, then $\widehat{\pi}_{w} \stackrel{a.s.}{\to} \pi_w $ when $n\to\infty$.
\end{theorem}

From the practical point of view, it is (weak) consistency that really matters, which also naturally leads to the exploration of asymptotic normality, and our following research path is in this direction. It leads us to practically attractive and justifiable conditions on the weight function $w$ as well as on the cdf's of $X$ and $Y$. Our focus now is also shifting from the simpler $\widehat{\pi}_w$ toward the more complex weighted allocation rule $\widehat{\Pi}_{w}$. Not surprisingly, therefore, in what follows we employ the conditional expectation function
\[
g_{X\mid Y}(y)=\mathbf{E}[X \mid Y=y]
\]
defined on the support of $Y$, as well as the conditional variance function
\[
v^2_{X\mid Y}(y)=\textbf{Var}[X \mid Y=y].
\]

We note that the function $g_{X\mid Y}\circ F_Y^{-1}(t)$ is known in the literature as the quantile regression function of $X$ on $Y$, and it has prominently manifested in the literature (e.g., Rao and Zhao, 1995; Tse, 2009; Tse, 2015; and the references therein). The quantile conditional-variance function $v^2_{X\mid Y}\circ F_Y^{-1}(t)$ is also prominently featured in these works. All these functions play a pivotal role throughout the rest of the present paper.

\begin{theorem}\label{theorem-2}
Assume that $\mathbf{E}[X^2]$ is finite and the cdf $F_Y$ is continuous. If $v^2_{X\mid Y}\circ F_Y^{-1}\in L_p$ and $w^2\in L_{q}$
for some $p,q \in [1, \infty ]$ such that $p^{-1}+q^{-1}=1$, then
$\widehat{\Pi}_{w} \stackrel{\mathbf{P}}{\to} \Pi_w $ when $n\to\infty$.
\end{theorem}

To appreciate the theorem from the practical perspective, we look at several special cases. First, when $p=1$, we have $q=\infty $, which is to say that the weight function $w$ is bounded. This covers the weight function $w(t)=\mathbf{1}\{t>p \}$ for every $p\in (0,1)$, and also the weight function $w(t)=\nu(1-t)^{\nu-1}$ for every $\nu\ge 1$. Note also that the condition $v^2_{X\mid Y}\circ F_Y^{-1}\in L_1$ is equivalent to $\mathbf{E}[X^2]<\infty$, which we assume.

Second, when $p=\infty$, which implies $q=1$, the function $v^2_{X\mid Y}(y)$ must be bounded and the  function $w^2$ integrable on $(0,1)$, that is, $w\in L_2$. The weight function $w(t)=\mathbf{1}\{t>p \}$ is always such, whereas $w(t)=\nu(1-t)^{\nu-1}$ belongs to $L_2$  only when $\nu>1-1/2$. The latter restriction has appeared naturally in Jones and Zitikis (2003, 2007), Brahimi et al. (2011), and other insurance-related works dealing with the proportional hazards premium.

The next theorem, which is in the spirit of Theorem \ref{theorem-zwet-2} and uses the notations introduced before it, concludes our explorations of consistency.

\begin{theorem}\label{theorem-3}
Assume that $\mathbf{E}[X^2]$ is finite and the cdf $F_Y$ is continuous.  If there is $\varepsilon>0$ such that  $v^2_{X\mid Y}\circ F_Y^{-1}\mathbf{1}_{B_{i,\varepsilon}}\in L_{p_i}$ and $w \mathbf{1}_{A_i}\in L_{2q_i}$ for every $i=1,\dots ,j$, then $\widehat{\Pi}_{w} \stackrel{\mathbf{P}}{\to} \Pi_w $ when $n\to\infty$.
\end{theorem}

We now set out to establish asymptotic normality of the estimator $\widehat{\Pi}_{w}$. We show, in particular, that its asymptotic variance is
\[
\sigma^2=
\big (\sigma_1^2+\sigma^2_{2}\big )/ \bigg(\int_0^1 w(u)\mathrm{d}u \bigg)^2
\]
with the notations
\begin{equation}
\sigma_1^2= \int_0^1v^2_{X\mid Y}\circ F_Y^{-1}(t) \,w^2(t)\, \mathrm{d}t
\label{sigma_1}
\end{equation}
and
\begin{equation}
\sigma_2^2
= \int_0^1 \int_0^1 w(s)w(t)\big(\min\{s,t\} -st\big ) \,\mathrm{d} g_{X\mid Y}\circ F_Y^{-1}(s)\,\mathrm{d} g_{X\mid Y}\circ F_Y^{-1}(t).
\label{sigma_2}
\end{equation}
We note at the outset that in the theorem that follows, we impose conditions that assure the finiteness of $\sigma_1^2$ and $\sigma_2^2$. Note also that for the variance $\sigma_2^2$ to be well defined, we need to, and thus do, assume -- without explicitly saying this every time -- that the quantile-regression function $g_{X\mid Y}\circ F_Y^{-1}$ is of bounded variation on $[\varepsilon,1-\varepsilon]$ for every $0<\varepsilon <1$.
Following Shorack's~(1972) terminology, this means that $g_{X\mid Y}\circ F_Y^{-1}$ belongs to the class $\mathcal{L}$. In general, every function $h\in \mathcal{L}$ generates a Lebesgue-Stieltjes signed measure whose total variation we denote by $|h|$.

\begin{theorem}\label{as-norm}
Assume that $\mathbf{E}[X^2]$ is finite and the cdf $F_Y$ is continuous. Furthermore, let the weight function $w$ satisfy the following conditions:
\begin{enumerate}[{\rm (i)}]
\item  \label{ia}
$w$ is continuous on $(0,1)$ except possibly at a finite number of points $t_1<\dots<t_m$, and there is $r>1/2$ such that for every $\varepsilon>0$ there is a constant $c<\infty $ such that
\[
|w(u)-w(v)|\leq c |u-v|^{r}
\]
for all $u,v\in (t_{i-1},t_i)\cap (\varepsilon, 1-\varepsilon)$ and every $i=1,\dots , m+1$, where $t_0:=0$ and $t_{m+1}:=1$;
\item  \label{iia}
there is (small) $\varepsilon>0$ such that $w$ is differentiable on the set $ \Theta_{\varepsilon}:=(0,\varepsilon)\cup(1-\varepsilon,1)$, and there are
$\kappa_1,\,\kappa_2 \in [0,1)$ such that
\[
t(1-t)|w'(t)|,\,|w(t)|\leq c  t^{-\kappa_1/2}(1-t)^{-\kappa_2/2}
\]
for all $t\in \Theta_{\varepsilon}$.
\end{enumerate}
If the function $g_{X\mid Y}\circ F_Y^{-1}$ is continuous at every point $t_i$ of condition (\ref{ia}), and, for some $\delta>0$, the bound
\begin{equation}
\label{thm-5_cond_g}
v^2_{X\mid Y}\circ F_Y^{-1}(t)\leq c  t^{-1+\kappa_1+\delta}(1-t)^{-1+\kappa_2+\delta}
\end{equation}
holds for all $t\in \Theta_{\varepsilon}$ with the same $\kappa_1$ and $\kappa_2 $ as in condition (\ref{iia}), then
\begin{equation}
n^{1/2}(\widehat{\Pi}_{w}-\Pi_{w})\stackrel{law}{\to}\mathcal{N}(0,\sigma^2)
\label{thm-5c_pi_w}
\end{equation}
when $n\to \infty $.
\end{theorem}

In the special case $Y=X$, we have the following corollary to Theorem \ref{as-norm}.

\begin{corollary}\label{as-norm-cor}
Assume that $\mathbf{E}[X^2]$ is finite and the cdf $F_X$ is continuous. If the weight function $w$ satisfies the conditions of Theorem \ref{as-norm} and the quantile function $ F_X^{-1}$  is continuous at every point $t_i$ of condition (\ref{ia}), then
\begin{equation}
n^{1/2}(\widehat{\pi}_{w}-\pi_{w})\stackrel{law}{\to}\mathcal{N}(0,\sigma^2)
\label{thm-5c_pi_w-cor}
\end{equation}
when $n\to \infty $, where
\[
\sigma^2
= {1\over \big(\int_0^1 w(u)\mathrm{d}u \big)^2}
\int_0^1 \int_0^1 w(s)w(t)\big(\min\{s,t\} -st\big ) \,\mathrm{d} F_X^{-1}(s)\,\mathrm{d} F_X^{-1}(t).
\]
\end{corollary}

In a variety of forms, Corollary \ref{as-norm-cor} has frequently appeared in literature. Indeed, $\widehat{\pi}_{w}$ is an $L$-statistic and $\pi_{w}$ is its asymptotic mean, called $L$-functional. For details and references on the topic, we refer to the monographs by Helmers (1982), Serfling (1980), and Shorack~(2017).

It is clear from the above results that the tails of the weight function $w$ and the cdf's $F_X$ and $F_Y$ interact, and thus there is always a ballancing act to maintain: stronger conditions on $w$ lead to weaker conditions on the cdf's, and vice versa. There is, however, a possibility to weaken both sets of conditions at the same time, but this leads to drastically different results and hinges on other techniques of proof, as seen from the works of Necir and Meraghni (2009), and Necir et al. (2007), who tackle the proportional hazards transform; Necir et al. (2010), and Rassoul (2013), who tackle the tail conditional expectation; and Brahimi et al. (2012), who tackle the general distortion risk measure.

\section{Proofs}
\label{sect-3a}

\paragraph{Proof of Theorem~\ref{theorem-1}.}

We have $\widetilde{\Pi}_{w} \stackrel{a.s.}{\to} \Pi_w $ provided that, when $n \to \infty$,
\begin{equation}
\label{s1}
{1\over n}\sum_{k=1}^{n}X_k w\circ \widehat{F}_Y(Y_k)\stackrel{a.s.}{\to} \mathbf{E}[Xw\circ F_Y(Y)]
\end{equation}
and
\begin{equation}
\label{s2}
{1\over n}\sum_{k=1}^{n} w\circ \widehat{F}_Y(Y_k)\stackrel{a.s.}{\to} \mathbf{E}[w\circ F_Y(Y)].
\end{equation}
Statement~\eqref{s2} follows from statement~\eqref{s1} if we set $X_k$ to $1$. Hence, we only need to prove statement~\eqref{s1}. We  write
\begin{equation}
\label{s1_p}
{1\over n}\sum_{k=1}^{n}X_k w\circ \widehat{F}_Y(Y_k)={1\over n}\sum_{k=1}^{n}X_k w\circ F_Y(Y_k)+{1\over n}\sum_{k=1}^{n}X_k \big ( w\circ \widehat{F}_Y(Y_k)-w\circ F_Y(Y_k) \big ).
\end{equation}
The classical strong law of large numbers implies that $n^{-1}\sum_{k=1}^{n}X_k w\circ F_Y(Y_k)$ converges to $\mathbf{E}[Xw\circ F_Y(Y)]$ almost surely. Hence, we are left to prove that the second average on the right-hand side of equation (\ref{s1_p}) converges to $0$ almost surely. This we achieve by first estimating its absolute value by
\begin{equation}
\label{s1_p3}
\bigg ( {1\over n}\sum_{k=1}^{n}|X_k| \bigg ) \sup_{y\in \mathbb{R}}\big |w\circ \widehat{F}_Y(y)-w\circ F_Y(y)\big |.
\end{equation}
By the strong law of large numbers, $n^{-1}\sum_{k=1}^{n}|X_k|$ converges almost surely to the (finite) mean of $|X|$, and the supremum in~(\ref{s1_p3}) converges to zero  almost surely because of the classical Glivenko-Cantelli Theorem and the uniform continuity of $w$, which holds because $w$ is continuous on the compact interval $[0,1]$.  Hence, statement \eqref{s1} holds, and so does $\widetilde{\Pi}_{w} \stackrel{a.s.}{\to} \Pi_w $. Statement $\widetilde{\pi}_{w} \stackrel{a.s.}{\to} \pi_w $ follows as a special case when $X=Y$ and $X_k=Y_k$ for all $k=1, \dots , n$. This completes the proof Theorem~\ref{theorem-1}.

\paragraph{Proof of Theorem \ref{theorem-2}.}

The theorem follows from the statement
\begin{equation}
\label{thm2_p1}
\widehat{\Delta}_{w}\stackrel{\mathbf{P}}{\to}
\mathbf{E}[Xw\circ F_Y(Y)].
\end{equation}
To prove it, we write the decomposition
$\widehat{\Delta}_{w}=(T_{n,1}+T_{n,2})/n$, where
\[
T_{n,1}=\sum_{k=1}^{n}g_{X\mid Y}(Y_{k:n}) w_{k,n}
\]
and
\[
T_{n,2}= \sum_{k=1}^{n}\big (X_{[k:n]} -g_{X\mid Y}(Y_{k:n})\big ) w_{k,n} .
\]
The rest of the proof consists of two steps:
\begin{gather}
{1\over n}T_{n,1} \stackrel{\mathbf{P}}{\to} \mathbf{E}[Xw\circ F_Y(Y)],
\label{proof_t2-0}
\\
{1\over n}T_{n,2} \stackrel{\mathbf{P}}{\to} 0.
\label{thm2_p1b}
\end{gather}

In fact, statement (\ref{proof_t2-0}) holds with convergence in probability replaced by convergence almost surely. Indeed, the strong law of large numbers for $L$-statistics (van~Zwet,~1980; Corollary~2.1) implies
\begin{equation}
\label{proof_t2}
{1\over n}T_{n,1} \stackrel{a.s.}{\to} \int_0^1 g_{X\mid Y}\circ F_Y^{-1}w\mathrm{d}\lambda
\end{equation}
when $n \to \infty$. It remains to note that the integral on the right-hand side of statement (\ref{proof_t2}) is equal to $\mathbf{E}[Xw\circ F_Y(Y)]$. Hence, we are left to prove  statement~\eqref{thm2_p1b}, which means that, for every $\delta >0$, we need to show
\begin{equation}
\label{thm2_p2b}
\mathbf{P}\big(|T_{n,2}| >n\delta \big) \rightarrow 0
\end{equation}
when $n\to \infty $. Recall that, conditionally on $Y_{1:n}, \dots, Y_{n:n}$, the concomitants $X_{[1:n]}, \dots , X_{[n:n]}$ are independent (Bhattacharya, 1974; Lemma 1). Hence, with the help of Markov's inequality, we obtain
\begin{align}
\mathbf{P}\big(|T_{n,2}| >n\delta \big)
&=\mathbf{E}\bigg [\mathbf{P}\Big(\Big|\sum_{k=1}^{n}\big (X_{[k:n]} -g_{X\mid Y}(Y_{k:n})\big ) w_{k,n}\Big| >n\delta \mid Y_1,\dots,Y_n\Big) \bigg ]
\notag
\\
&\le \frac{1}{n^2\delta^2}\sum_{k=1}^{n}\mathbf{E}\Big [\mathbf{E} \Big [\big (X_{[k:n]} -g_{X\mid Y}(Y_{k:n})\big )^2\mid Y_1,\dots,Y_n\Big] \Big ] w^2_{k,n}
\notag
\\
&= \frac{1}{n^2\delta^2}\sum_{k=1}^{n}\mathbf{E}\big [ v^2_{X\mid Y}(Y_{k:n})\big] w^2_{k,n} ,
\label{thm2_p3b}
\end{align}
where the right-most equation follows from the fact that (Bhattacharya, 1974; Lemma 1) conditionally on $Y_{1:n}, \dots, Y_{n:n}$, the concomitants $X_{[1:n]}, \dots , X_{[n:n]}$ follow the cdf's $F(x\mid Y_{1:n}), \dots , F(x\mid Y_{n:n})$, respectively, where
$F(x\mid y)=\textbf{P}[X\leq x\mid Y=y]$. Next we apply H\"{o}lder's inequality on the right hand side of bound \eqref{thm2_p3b} and obtain
\begin{align}
\mathbf{P}\big(|T_{n,2}| >n\delta \big)
&\le \frac{1}{n\delta^2}
\bigg(\frac{1}{n}\sum_{k=1}^{n}\Big ( \mathbf{E}\big [ v^2_{X\mid Y}(Y_{k:n})\big]\Big )^p\bigg)^{1/p}  \bigg({1\over n} \sum_{k=1}^{n} |w_{k,n}|^{2q}\bigg)^{1/q}
\notag
\\
&\le  \frac{1}{n\delta^2}
\bigg(\frac{1}{n}\sum_{k=1}^{n}\mathbf{E}\big [ v^{2p}_{X\mid Y}(Y_{k:n})\big]\bigg)^{1/p}  \bigg({1\over n} \sum_{k=1}^{n} |w_{k,n}|^{2q}\bigg)^{1/q}
\notag
\\
&=  \frac{1}{n\delta^2}
\bigg(\frac{1}{n}\sum_{k=1}^{n}\mathbf{E}\big [ v^{2p}_{X\mid Y}\circ F^{-1}_Y (U_{k})\big]\bigg)^{1/p}  \bigg({1\over n} \sum_{k=1}^{n} |w_{k,n}|^{2q}\bigg)^{1/q},
\label{thm2_p5b}
\end{align}
where $U_1,\dots,U_n$ are independent and $(0,1)$-uniform random variables.
By the classical law of large number, the first sum on the right-hand side of equation (\ref{thm2_p5b}) convergence to the integral $\int_0^1 v^{2p}_{X\mid Y}\circ F^{-1}_Y (t)\,\mathrm{d}t$, whereas the second sum converges to $\int_0^1 |w(t)|^{2q}\mathrm{d}t$. Both integrals are finite by assumption. This completes the proof of Theorem \ref{theorem-2}.

\paragraph{Proof of Theorem \ref{theorem-3}.}

We need to prove statement~\eqref{thm2_p1} under the conditions of Theorem~\ref{theorem-3}. We  start again with the decomposition $\widehat{\Delta}_{w}=(T_{n,2}+T_{n,1})/n$. Statement~\eqref{proof_t2} follows by the strong law of large numbers for $L$-statistics (van~Zwet, 1980; Theorem~3.1). It remains to prove statement \eqref{thm2_p1b}. We fix any $\delta>0$ and write
\begin{align}
\mathbf{P}\big(|T_{n,2}| >n\delta \big)
=&\mathbf{P}\bigg(\Big|\sum_{k=1}^{n}\big (X_{[k:n]} -g_{X\mid Y}\circ F^{-1}(U_{k:n})\big ) w_{k,n}\Big | >n\delta \bigg)
\notag
\\
\leq &\Delta+\mathbf{P}\big(\mathcal{D}^c\big)
\label{thm3_p3}
\end{align}
for any subset $\mathcal{D}$ of the sample space, where
\[
\Delta:=\mathbf{P}\bigg(\bigg\{\Big|\sum_{k=1}^{n}\big (X_{[k:n]} -g_{X\mid Y}\circ F^{-1}(U_{k:n})\big ) w_{k,n}\Big | >n\delta\bigg\}\cap \mathcal{D} \bigg)
\]
and $U_{1:n}, \dots , U_{n:n}$ are the order statistics corresponding to  independent $(0,1)$-uniform random variables $U_1,\dots,U_n$. We next make a special choice of $\mathcal{D}$.

First, we recall the definitions of $A_i$ and $B_{i,\varepsilon}$ that are given before Theorem \ref{theorem-zwet-2}. Then we define
$r_i=\min\{k:\, k/n\in A_i\}$ and $ s_i=\max\{k:\, k/n\in A_i\}$, and with the notation $\mathcal{D}_i= \big\{ U_{r_i:n}\in B_{i,\varepsilon}\big\}\cap \big\{ U_{s_i:n}\in B_{i,\varepsilon}\big\}$, we define
\[
\mathcal{D}=\bigcap_{i=1}^{j}\mathcal{D}_i.
\]
Since $r_i/n,\, s_i/n\in A_i \subset B_{i,\varepsilon}$,  Bernstein's inequality implies $\mathbf{P}\big(\mathcal{D}^c_i\big)\leq \exp\{-c_i n\}$ for some $c_i>0$, where
\[
\mathcal{D}^c_i=\big\{ U_{r_i:n}\notin B_{i,\varepsilon}\big\}\cup \big\{ U_{s_i:n}\notin B_{i,\varepsilon}\big\}.
\]
Consequently,
\begin{equation}
\label{thm3_p2}
\mathbf{P}\big(\mathcal{D}^c\big) =\mathbf{P}\bigg (\bigcup_{i=1}^{j}\mathcal{D}_i^c \bigg) \leq m\exp\{-cn\}
\end{equation}
with $c=\min_{1\leq i\leq j}c_i$. In view of estimate~\eqref{thm3_p2}, from now on we restrict our attention to only the quantity $\Delta$.

Since conditionally on $Y_{1:n}, \dots, Y_{n:n}$, the concomitants $X_{[1:n]}, \dots , X_{[n:n]}$ follow the cdf's $F(x\mid Y_{1:n}), \dots , F(x\mid Y_{n:n})$, respectively, we use Markov's inequality and obtain the bound
\begin{align}
\Delta
&=\mathbf{E}\bigg [ \mathbf{1}_{\mathcal{D}}\mathbf{P}\bigg(\Big |\sum_{k=1}^{n}\big (X_{[k:n]} -g_{X\mid Y}\circ F^{-1}(U_{k:n})\big ) w_{k,n}\Big| >n\delta \mid U_1,\dots,U_n \bigg)\bigg ]
\notag
\\
&\leq \frac{1}{(n\delta)^2}\mathbf{E}\bigg [\mathbf{1}_{\mathcal{D}}\sum_{k=1}^n v^2_{X\mid Y}\circ F^{-1}_Y (U_{k:n}) w^2_{k,n} \bigg ].
\label{thm3_p4-i}
\end{align}

We split the sum $\sum_{k=1}^n$ into $\sum_{i=1}^{j}  \sum_{k:k/n=a_i}$ and $ \sum_{i=1}^j \sum_{k:\, k/n\in A_i}$. (If the sum $\sum_{k:k/n=a_i} $ is empty, we set it to $0$ by definition.) Hence, in order to show that $\Delta$ converges to $0$ when $n\to \infty $, we need to show that, for every $i=1,\dots , j$, the right-hand side of equation (\ref{thm3_p4-i}) converges to $0$ when the sum $\sum_{k=1}^n$ is replaced by $\sum_{k:k/n=a_i}$ as well as by $\sum_{k:\, k/n\in A_i}$. We begin with the first sum $\sum_{k:k/n=a_i}$, which is either empty or contains only one summand. If it is not empty, then let $k$ be the (only) integer that satisfies $k=na_i$, thus obviously implying $na_i\in \mathbb{N}$. In this case, we have
\begin{align}
\frac{1}{n}\mathbf{E}\bigg [\mathbf{1}_{\mathcal{D}}
& v^2_{X\mid Y}\circ F^{-1}_Y (U_{k:n}) w^2_{k,n} \bigg ]
\notag
\\
&\le
\frac{c}{n}\mathbf{E}\big [v^2_{X\mid Y}\circ F^{-1}_Y (U_{na_i:n})\big ]
\notag
\\
&= \frac{c}{n}\int_0^1 v^2_{X\mid Y}\circ F^{-1}_Y (t)  \frac{n!}{(na_i-1)!(n-na_i)!}t^{na_i-1}(1-t)^{n-na_i} \, \mathrm{d}t
\notag
\\
&\le c \int_0^1 v^{2}_{X\mid Y}\circ F^{-1}_Y (t) \sum_{k=1}^{n}\frac{(n-1)!}{(k-1)!(n-k)!}t^{k-1}(1-t)^{n-k} \, \mathrm{d}t
\notag
\\
&= c \int_0^1 v^{2}_{X\mid Y}\circ F^{-1}_Y (t)  \mathrm{d}t <\infty .
\label{thm3_p4-a}
\end{align}
Consequently, we are now left to deal only with the sum $\sum_{k:\, k/n\in A_i}$. That is, we conclude that $\Delta$ converges to $0$ if, for every $i=1,\dots , j$, the quantity
\[
\Delta_{2,i}
:=\frac{1}{(n\delta)^2} \sum_{k:\, k/n\in A_i} \mathbf{E}\big [\mathbf{1}_{\mathcal{D}}v^2_{X\mid Y}\circ F^{-1}_Y (U_{k:n})\big ] w^2_{k,n}
\]
converges to $0$ when $n\to \infty $. With the help of H\"{o}lder's inequality, we have
\begin{multline}
\Delta_{2,i}
\leq \frac{1}{(n\delta)^2}
\bigg(\sum_{k:\, k/n\in A_i}\Big ( \mathbf{E}\big [\mathbf{1}_{\mathcal{D}}v^2_{X\mid Y}\circ F^{-1}_Y (U_{k:n})\big ]\Big )^{p_i}\bigg)^{1/p_i}
\\
\times \bigg(\sum_{k:\, k/n\in A_i} |w_{k,n}|^{2q_i}\bigg)^{1/q_i}.
\label{thm3_p4}
\end{multline}
We have, when $n\to \infty $,
\[
{1\over n} \sum_{k:\, k/n\in A_i} |w_{k,n}|^{2p_i}\to \int_{A_i} |w(t)|^{2p_i}\, \mathrm{d}t <\infty .
\]
Furthermore,
\begin{align}
\frac{1}{n}\sum_{k:\, k/n\in A_i}&\Big ( \mathbf{E}\big [\mathbf{1}_{\mathcal{D}}v^2_{X\mid Y}\circ F^{-1}_Y (U_{k:n})\big ]\Big )^{p_i}
\notag
\\
&\le \frac{1}{n}\sum_{k:\, k/n\in A_i}\mathbf{E}\big [\mathbf{1}_{\mathcal{D}_i}v^{2p_i}_{X\mid Y}\circ F^{-1}_Y (U_{k:n})\big ],
\notag
\\
&= \sum_{k:\, k/n\in A_i}\int_{B_{i,\varepsilon}}v^{2p_i}_{X\mid Y}\circ F^{-1}_Y (t) \frac{(n-1)!}{(k-1)!(n-k)!}t^{k-1}(1-t)^{n-k} \, \mathrm{d}t
\notag
\\
&\le \int_{B_{i,\varepsilon}}v^{2p_i}_{X\mid Y}\circ F^{-1}_Y (t) \sum_{k=1}^n \frac{(n-1)!}{(k-1)!(n-k)!}t^{k-1}(1-t)^{n-k} \, \mathrm{d}t
\notag
\\
&= \int_{B_{i,\varepsilon}}v^{2p_i}_{X\mid Y}\circ F^{-1}_Y (t) \, \mathrm{d}t<\infty.
\label{bound-1}
\end{align}
Bounds~\eqref{thm3_p3}--\eqref{bound-1} imply statement \eqref{thm2_p1b}, thus completing the proof of Theorem~\ref{theorem-3}.

\paragraph{Proof of Theorem \ref{as-norm}.}

Asymptotic normality of $\widehat{\Pi}_{w}$ follows if we show that, when $n\to \infty $,
\begin{gather}
n^{1/2}(\widehat{\Pi}_{w}-\Pi_{w,n})\stackrel{law}{\to}\mathcal{N}(0,\sigma^2),
\label{thm-5a_pi_w}
\\
n^{1/2} (\Pi_{w,n}-\Pi_w)\to 0,
\label{thm-5b_pi_w}
\end{gather}
where
\[
\Pi_{w,n}:=\sum_{k=1}^n {w_{k,n}\over \int_0^1 w(u)\mathrm{d}u} \int_{(k-1)/n}^{k/n} g_{X\mid Y}\circ F_Y^{-1}(t)\, \mathrm{d}t.
\]
We next establish statements (\ref{thm-5a_pi_w}) and (\ref{thm-5b_pi_w}), and in this way complete the proof of Theorem \ref{as-norm}. We note that the two statements require different subsets of conditions formulated in Theorem \ref{as-norm}. In the proofs that follow, we shall specify which of them, and where, are required.

\paragraph{Proof of statement (\ref{thm-5a_pi_w}).}

We have $\Pi_{w,n}=\Delta_{w,n}/\int_0^1 w(u)\mathrm{d}u$, where
\[
\Delta_{w,n}=\sum_{k=1}^n w_{k,n} \int_{(k-1)/n}^{k/n} g_{X\mid Y}\circ F_Y^{-1}(t)\, \mathrm{d}t
\]
Hence, statement (\ref{thm-5a_pi_w}) is equivalent to
\[
n^{1/2}(\widehat{\Delta}_{w}-\Delta_{w,n})
\stackrel{law}{\to}\mathcal{N}(0,\sigma_1^2+\sigma^2_{2})
\]
when $n\to \infty $. We write
$n^{1/2}(\widehat{\Delta}_{w}-\Delta_{w,n})=W_n + T_n$,
where
\begin{equation}
\label{thm_4_p_1}
W_n=n^{-1/2} \sum_{k=1}^n \big (X_{[k:n]} -g_{X\mid Y}(Y_{k:n})\big ) w_{k,n}
\end{equation}
and
\begin{equation}
\label{thm_4_p_2}
T_n= n^{1/2} \bigg({1\over n}\sum_{k=1}^n g_{X\mid Y}(Y_{k:n}) w_{k,n} -\Pi_{w,n}\bigg).
 \end{equation}
Hence, to prove the theorem, we need to show that
\begin{equation}
W_n +T_n\stackrel{law}{\to}\mathcal{N}(0,\sigma_1^2+\sigma^2_{2}),
\label{thm-5f_pi_w}
\end{equation}
when $n\to \infty $. We follow the approach of Yang~(1981) for proving the central limit theorem for linear combinations of concomitants.


\begin{theorem}[Yang, 1981]\label{lem_yang}
Let $(X_1,Y_1),(X_2,Y_2),\dots$ be random pairs and, for every
$n\geq 1$, the first $n$ pairs $(X_1,Y_1),\dots,(X_n,Y_n)$
possess a joint distribution. Denote  $\mathbf{Z}_n= ((X_1,Y_1),\dots,(X_n,Y_n))$ and
 $\mathbf{Y}_n= (Y_1,\dots,Y_n)$, and let $W_n:={W}_n(\mathbf{Z}_n)$ and $T_n:={T}_n(\mathbf{Y}_n)$ be measurable vector-valued functions of $\mathbf{Z}_n$ and $\mathbf{Y}_n$,  respectively. Suppose ${T}_n$ converges in distribution to $F_T$, and the conditional distribution of ${W}_n$ given $\mathbf{Y}_n$ converges weakly to a distribution $F_W$ which does not depends on the $Y_k$'s. Then
 $({W}_n,{T}_n)\stackrel{law}{\rightarrow}F_W F_T$.
\end{theorem}

First we work with the quantity $W_n$ defined by equation \eqref{thm_4_p_1}, and  prove that its conditional distribution given $\mathbf{Y}_n $ tends to the normal distribution with the mean $0$ and variance   $\sigma_1^2$ for almost all sequences  $(Y _m)_{m\ge 1}$, with the limiting distribution not depending on the sequence $(Y_m)_{m\ge 1}$. Next, we prove that the quantity $T_n$  defined by equation \eqref{thm_4_p_2} is asymptotically normal with the mean $0$ and variance $\sigma^2_{2}$. Given these two results, Theorem \ref{lem_yang} implies that the joint distribution of $(W_n,T_n)$ converges to the product of the two aforementioned normal distributions. In turn, this implies that $W_n+T_n$ is asymptotically normal with the mean $0$ and variance  $\sigma_1^2+\sigma^2_{2}$. Hence, the rest of the proof consists of two parts, and they deal with the asymptotic normality of $W_n$ and $T_n$, respectively.

\paragraph{\it Part 1.}

Using Bhattacharya's (1974) result already utilized in the proof of Theorem~\ref{theorem-3}, we have
$\mathbf{E}[W_n \mid \boldsymbol{Y}_n]=0$ with the conditional variance $V_n^2:=\textbf{Var}[W_n\mid \boldsymbol{Y}_n]$ expressed by
\[
V_n^2={1\over n}\sum_{k=1}^{n}v^2_{X\mid Y}(Y_{k:n}) w^2_{k,n}.
\]
Applying Lindeberg's normal-convergence criterion, we conclude that the sequence of the (conditional) distributions of $W_n/V_n$ is asymptotically standard normal if, for every $\varepsilon>0$ and when $n\to \infty $,
\begin{equation}
\label{thm_4_p_3}
\frac 1{nV_n^2}\sum_{k=1}^{n} w^2_{k,n}  h_{\theta_{k,n}}(Y _{k:n}) \rightarrow 0
\end{equation}
for almost all realizations of the sequence $Y_1, Y_2, \dots $, where
\begin{equation}
\label{thm_4_p_4}
h_{\theta_{k,n}}(y)=\int (x-g_{X \mid Y}(y))^2 \mathbf{1}\{|x-g_{X \mid Y}(y)|\geq \theta_{k,n}\} \mathrm{d}F(x\mid y)
\end{equation}
with the notation
\[
\theta_{k,n}={\varepsilon n^{1/2} V_n \over \left|w_{k,n}\right|}.
\]
(If $w_{k,n}=0$, the  corresponding summand in statement \eqref{thm_4_p_3} vanishes, and hence $\theta_{k,n}$ can be defined arbitrarily in this case.) The strong law of large numbers for $L$-statistics (van~Zwet, 1980; Theorem~3.1) implies
\begin{equation}
\label{thm_4_p_5}
V_n^2  \stackrel{a.s.}{\to} \int_0^1v^2_{X \mid Y}\circ  F_Y^{-1}(t) w^2(t)\, \mathrm{d}t ,
\end{equation}
with the integral on the right-hand side being equal to $\sigma_1^2 $. To verify $\theta_{k,n}\stackrel{a.s.}{\to} \infty$, we write the bounds
\begin{align}
\theta_{k,n}
&\geq \varepsilon n^{1/2} V_n / \max_{k=1,\dots,n}|w_{k,n}|
\notag
\\
&\geq \varepsilon n^{1/2} V_n /\max_{k=1,\dots,n}
\bigg ( {k(n-k)\over n^2} \bigg )^{-\max(\kappa_1,\kappa_2)/2}
\notag
\\
&= \varepsilon n^{1/2} V_n / n^{\max(\kappa_1,\kappa_2)/2}.
\label{thm_4_p_6}
\end{align}
Since $\max(\kappa_1,\kappa_2)<1$, we have $\theta_{k,n}\stackrel{a.s.}{\to} \infty$. Applying the strong law of large numbers for $L$-statistics (van~Zwet, 1980; Theorem~3.1), we have that, for every $K>0$,
\begin{equation}
\label{thm_4b_p_5}
\frac 1{n}\sum_{k=1}^{n} w^2_{k,n} h_{K}(Y_{k:n}) \stackrel{a.s.}{\to}
\int_0^1  w^2(t)\, h_K\circ F_Y^{-1}(t)\, \mathrm{d}t
\end{equation}
when $n \to \infty$. (The function $h_K(y)$ is defined by equation \eqref{thm_4_p_4} with $K$ instead of $\theta_{k,n}$.) Since  $\theta_{k,n}\to \infty$, statements (\ref{thm_4b_p_5}) and \eqref{thm_4_p_5} imply  the Lindeberg's criterion for almost all realizations of the sequence $(Y_m)_{m\ge 1}$. Hence, the conditional distribution of $W_n/V_n$ given $\boldsymbol{Y}_n$ converges to the standard normal distribution almost surely.

\paragraph{\it Part 2.}
In order to prove statement (\ref{thm-5a_pi_w}), it remains to show that the distribution of $T_n$ given by~\eqref{thm_4_p_2} tends to the normal law with the mean $0$ and variance $\sigma^2_2$. The latter fact is a direct convergence of  a result of Shorack's (1972) Theorem~1 on asymptotic normality of the linear combination of functions of order statistics.

Indeed, let $U_1,U_2,\dots$ be a sequence of independent $(0,1)$-uniformly distributed random variables, and let $U_{k:n}$, $k=1,\dots,n$, denote the order statistics based on the first $n$ members of the sequence. Then, with the equality holding in distribution, we have
\[
T_n=n^{1/2} \bigg({1\over n}\sum_{k=1}^n g_{X\mid Y}\circ F_Y^{-1}(U_{k:n}) w_{k,n} -\Delta_{w,n}\bigg).
\]
Since under the conditions of  Theorem \ref{as-norm}, the conditions  of Theorem~1 by Shorack~(1972) are satisfied, the aforementioned asymptotic normality of $T_n$ holds.  Statement (\ref{thm-5a_pi_w}) follows.

\paragraph{Proof of statement (\ref{thm-5b_pi_w}).}
We start with the equations
\begin{align}
\Delta_{w,n}&=\sum_{k=1}^n w_{k,n} \int_{(k-1)/n}^{k/n} g_{X\mid Y}\circ F_Y^{-1}(t)\, \mathrm{d}t
\notag
\\
&=\int_0^1 g_{X\mid Y}\circ F_Y^{-1}(t) w_n(t)\, \mathrm{d}t
\label{thm-4_01}
\end{align}
where the function $w_n:\ (0,1]\to \mathbb{R}$ is defined by
$w_n(t)=w_{k,n}$ when $(k-1)/ n <t\le k/n $, for all $k=1,\dots , n$. Next we write
\begin{equation*}
\label{thm_5_p_1}
n^{1/2} (\Pi_{w,n}-\Pi_w)=I_{n,1}+I_{n,2}+I_{n,3},
\end{equation*}
where
\[
I_{n,l}=n^{1/2}\int_{D_l} \big(w_n(t)-w(t)\big)g_{X\mid Y}\circ F_Y^{-1}(t)\, \mathrm{d}t
\]
with the sets
$D_1=(0,\varepsilon)$,  $D_2=(\varepsilon, 1-\varepsilon)$, and $D_3=(1-\varepsilon,1)$, and with a sufficiently small $\varepsilon>0$ so that we could use condition (\ref{iia}) of Theorem \ref{as-norm}.
We shall prove that $I_{n,1}$ and $I_{n,2}$ converge to zero when $n\to \infty$; the treatment of $I_{n,3}$ is similar as in the case of $I_{n,1}$ and is therefore omitted. We have
\begin{multline}
I_{n,1}=n^{1/2} \int_0^{1/n} \big(w_n(t)-w(t)\big)g_{X\mid Y}\circ F_Y^{-1}(t)\, \mathrm{d}t
\\
+ n^{1/2}\int_{1/n}^{\varepsilon} \big(w_n(t)-w(t)\big)g_{X\mid Y}\circ F_Y^{-1}(t)\, \mathrm{d}t.
\label{thm_5_p_2}
\end{multline}

The first integral on the right-hand side of equation \eqref{thm_5_p_2} is equal to
\begin{equation}
\label{thm_abc0}
w(1/n) \int_0^{1/n}g_{X\mid Y}\circ F_Y^{-1}(t)\, \mathrm{d}t -  \int_0^{1/n} w(t)g_{X\mid Y}\circ F_Y^{-1}(t)\, \mathrm{d}t.
\end{equation}
By condition (\ref{iia}) of Theorem \ref{as-norm}, the absolute value of quantity (\ref{thm_abc0}) does not exceed $c\, n^{-1/2-\delta/2}$.

Up to a constant, the absolute value of the second integral on the right-hand side of equation \eqref{thm_5_p_2} does not exceed
\[
\Theta_n:={1\over n} \sum_{k=2}^{[n\varepsilon]+1}  \int_{(k-1)/n}^{k/n}\tau_{k}^{-\kappa_1/2-1}\,t^{-1/2+\kappa_1/2+\delta/2}\, \mathrm{d}t
\]
for some $\tau_{k}\in \big((k-1)/n,\,k/n\big)$, where $[\cdot ]$ denotes the integer part. Without loss of generality, we let $\delta>0$ be smaller than $1-\kappa_1$. Next, we estimate $\Theta_n$ as follows:
\begin{align*}
\Theta_n
&\leq {c\over n}\Big( n^{1/2-\delta/2} + {1\over n}\sum_{k=2}^{[n\varepsilon]}\Big(\frac{k}{n}\Big)^{-3/2+\delta/2} \Big)
\\
&\leq  {c\over n}\Big( n^{1/2-\delta/2} + \int_{1/n}^{\varepsilon}\, t^{-3/2+\delta/2}\, \mathrm{d}t \Big)
\\
&\leq c\, n^{-1/2-\delta/2}.
\end{align*}
This bound yields $I_{n,1}\to 0$ when $n\to\infty$. It remains to prove that  $I_{n,2}$ converges to $0$ when $n\to \infty $. For this, we first rewrite $I_{n,2}$ as follows
\begin{equation*}
\label{thm_5_p_5}
I_{n,2}=n^{1/2} \sum_{i=1}^{m+1}  \int_{t_{i-1}}^{t_i} \big(w_n(t)-w(t)\big)g_{X\mid Y}\circ F_Y^{-1}(t)\, \mathrm{d}t.
\end{equation*}
By condition (\ref{ia}) of Theorem \ref{as-norm} and the integrability of $g_{X\mid Y}\circ F_Y^{-1}$ (which holds because the first moment of $X$ is finite), the  absolute value of $I_{n,2}$ does not exceed $c n^{1/2-r} $, which tends to $0$ because of~$r>1/2$. This completes the proof  statement (\ref{thm-5b_pi_w}).

Having thus established statements (\ref{thm-5a_pi_w}) and (\ref{thm-5b_pi_w}), we conclude the proof of Theorem~\ref{as-norm}.




\end{document}